\documentclass{article}%
\usepackage{amsmath, amsfonts,amssymb,graphicx}%
\usepackage{multicol}

\newtheorem{thm}{Theorem}

\newtheorem{pro}[thm]{Proposition}
\newtheorem{con}[thm]{Conjecture}
\newtheorem{rem}[thm]{Remark}

\newtheorem{que}[thm]{Question}

\def\Bbb{\mathbb}

\def\a{\alpha}
\def\b{\beta}

\def\per{\mathop{\rm per\,}}
\def\1per{\mathop{\mbox{\rm 1-per\,}}}
\def\2per{\mathop{\mbox{\rm 2-per\,}}}
\def\kper{\mathop{\mbox{\rm $k$-per\,}}}
\def\Per{\mathop{\rm Per\,}}

\begin{document}
\title{\Large The permanent functions of tensors}
\author{Qing-wen Wang${}^{\,\rm a}$%
,\;
Fuzhen Zhang${}^{\,\rm b}$
\\
\footnotesize{${}^{\,\rm a}$ Shanghai University, Shanghai, P.R. China; wqw@t.shu.edu.cn}\\
\footnotesize{${}^{\,\rm b}$ Nova Southeastern University, Fort Lauderdale, USA; zhang@nova.edu}}
\date{}
\maketitle

\bigskip
 \hrule
\bigskip

\noindent {\bf Abstract}
By a tensor  we mean a multidimensional array (matrix) or hypermatrix over a number field. This article aims to set an account of the studies on the permanent functions of tensors. We formulate the definitions of 1-permanent, 2-permanent, and $k$-permanent of a tensor in terms of hyperplanes, planes and $k$-planes of the tensor; we discuss the polytopes of stochastic tensors; at end  we present an extension of the generalized matrix function for tensors.

\medskip
{\footnotesize

}

\bigskip
\hrule


\section{Introduction}

 The study on
multidimensional arrays (or matrices)  may date back as early as the nineteenth century by Cayley
\cite{Cay1845,Cay1849}. Jurkat and Ryser revived  the topic in their seminal paper \cite{{JurRys68}} in 1968 in which they investigated configurations
and decompositions for multidimensional arrays.  Jurkat and Ryser's work was followed by a great deal of research on the topic, mainly on the combinatorial aspects of  certain types (such as stochasticity)  of multidimensional arrays; see, e.g., Brualdi and Csima \cite{BrCs75, BrCs91}.
In recent years, multidimensional arrays are found   applications in practical fields such as image processing (see, e.g.,  Qi and Luo \cite{Qi2017Book}), theory of computing (see, e.g.,
Cifuentes and Parrilo \cite{CiPa16}), and physics (see, e.g., Tichy \cite{Tic15}).
We are concerned with the  permanent functions of multidimensional arrays. Our purpose is to set an account on the specific topic based on  publications, including, in particular,  the ones
by Dow and Gibson \cite{DowGib87} and Taranenko \cite{Tar16}. The results are  expositorily presented with explanations other than in the format of theorem-proofs. Some results are easy observations; they are not necessarily new.
For the determinants of multidimensional arrays, hyperdeterminants, and related topics,
see, e.g.,  \cite{Gel94, Hu2013, Lim13, SSZ13}.

Let ${n_1, n_2, \dots, n_d}$ be positive integers. We write $A=(a_{i_1i_2\dots i_d})$,
$i_k=1, 2,  \dots, n_k$, $k=1, 2,  \dots, d$,  for an $n_1\times n_2\times  \cdots \times n_d$  multidimensional array or hypermatrix  of order $d$ (the number of indices).
Multidimensional arrays, or hypermatrices, or matrices of higher orders, are
referred to as {\em tensors}; see, e.g., \cite{DW16, KB09,  Qi2017Book}.  So, by  a tensor
  we mean a multidimensional array.
The tensors of order 1 (i.e., $d=1$) are vectors in $\Bbb R^{n_1}$, while the 2nd order tensors are just regular $n_1\times n_2$ matrices. A 3rd order tensor, i.e., an $n_1\times n_2\times n_3$ tensor,  may be viewed as a {book} of  $n_3$ pages (slices), each page is
an $n_1\times n_2$ matrix.

If $n_1= n_2=\dots= n_d=n$, we say that
$A$ is of order $d$ and dimension $n$  or we say that $A$ is an $\overbrace{n\times   \cdots \times n}^d$ tensor.
We also call an $n\times n\times n$ tensor (i.e., of order 3 and dimension $n$)
a {\em tensor cube} or
a {\em 3D matrix}.
We  refer to the  permanents of multidimensional arrays as the  {\em permanents} of
tensors, or hyperpermanents. Following the line of Dow and Gibson \cite{DowGib87}, we will begin with the definitions of 1-permanent, 2-permanent, and  $k$-permanent of tensors. 1-permanents, the most modest ones,  are useful in studying  hypergraphs
(see, e.g.,
 \cite{DowGib87, Tar16}), 
the 2-permanents with $d=3$ or of special relation of $d$ and $n$  are found applications in projective planes (see, e.g.,
 \cite{DowGib87}) 
  and polytope theory (see, e.g., \cite{ChangPZ16, CLN14, LZZ17}),
  while $k$-permanents are certainly an object in combinatorics themselves.

  \begin{rem}{\rm
We adopt Lim's terminology in \cite{Lim13} (see also \cite{KB09, Qi2017Book}), calling an $\overbrace{n\times   \cdots \times n}^d$ tensor a tensor of {\em order} $d$ and  {\em dimension} $n$.
 Such a tensor is also said
to be of order $n$ and  dimension $d$ in the literature; see, e.g., \cite{BrCs91, Tar16}.
  }
  \end{rem}

\section{The definitions of permanents of tensors}

\subsection{1-permanent and 2-permanent}
Let
$A=(a_{i_1\dots i_d})$ be an $n_1\times \cdots \times n_d$ tensor of order $d$ with real entries.
Dow and Gibson   \cite{DowGib87} defined (over a commutative ring) the
{\em permanent} of $A$   as
\begin{equation}\label{Def:1-per}
\per (A)=\sum \prod_{i=1}^{n_1} a_{i\sigma_2(i)\dots \sigma_d(i)},
\end{equation}
where the summation runs over all one-to-one functions $\sigma_k$ from
$\{1, 2, \dots, n_1\}$ to $\{1, 2, \dots, n_k\}$, $k=2, 3, \dots, d$, with $\per (A)=0$ if
$n_1>n_k$ for some $k$.

Note: under the definition (\ref{Def:1-per}),  if $A$ is an $n_1\times n_2$ matrix  and  $n_1>n_2$, then $\per (A)=0$, but $\per (A^t)$ need not be   0, where $A^t$ is the transpose of $A$. This is not in agreement with the fact that a matrix (i.e., order 2 tensor) and its transpose have the same permanent.
We may slightly modify and extend the definition (\ref{Def:1-per}) as follows. Let $n=\min\{n_1, n_2, \dots, n_d\}=n_j$ for some $j$. Then
\begin{equation}\label{Def:1-perA}
\per (A)=\sum \prod_{i=1}^{n} a_{\sigma_1(i)\cdots
\sigma_j(i)\dots \sigma_d(i)},
\end{equation}
where the summation runs over all one-to-one functions $\sigma_k$ from
$\{1, 2, \dots, n\}$ to $\{1, 2, \dots, n_k\}$, $k\not =j$, and $\sigma_j$ is the identity map. (\ref{Def:1-perA}) reduces to (\ref{Def:1-per}) if $n=n_1$. It is immediate by definition (\ref{Def:1-perA})  that $\per (A)=\per (A^t)$ for rectangular matrices.

When $n_1=n_2=\cdots =n_d=n$, (\ref{Def:1-per}) can be written in a symmetric form
\begin{equation}\label{Def:1-perB}
\per (A)=\frac{1}{n!}\sum_{\pi_1, \dots, \pi_d\in S_n} \prod_{i=1}^n a_{\pi_1(i)\cdots  \pi_d(i)},
\end{equation}
where $S_n$ is the symmetric group of degree $n$.

If $d=2$,   (\ref{Def:1-perB})  reduces to the usual permanent for
square matrices.

The definition (\ref{Def:1-per}) of permanent is in fact the so-called  {\em 1-permanent} (or 1-per for short) of the tensor $A$. 1-per (v.s. $k$-per; see the definition below or see
\cite[Sec.\,4]{DowGib87}) of $A$ is the sum of all  products of   $n_1$ entries of $A$ no two of which are taken from the same {\em hyperplane} (of order $d-1$; see Sec.\,2.2).  In the case of 3D matrices,
the {\em planes} of $A$ are the submatrices obtained by fixing one of $i,j,k$, and
the {\em lines} of $A$ are the submatrices obtained by fixing two of $i,j,k$.
\begin{figure*}[h]
\begin{center}
\begin{multicols}{2}
   \includegraphics[width=1.6in,totalheight=1.3in]{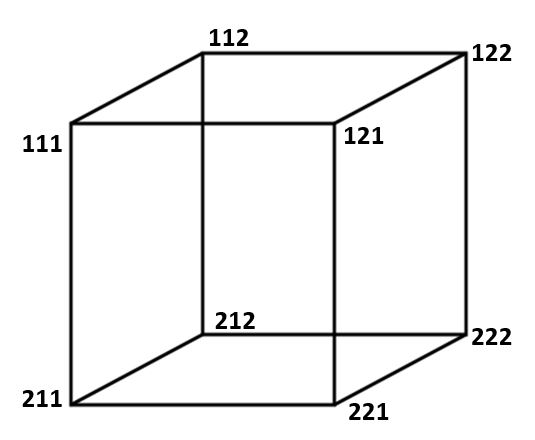}\par
    \includegraphics[width=1.5in,totalheight=1.2in]{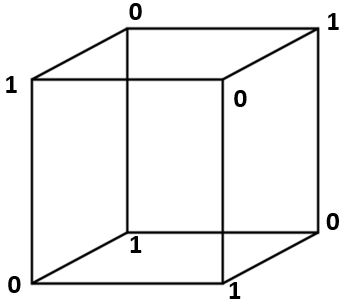}\par
    \end{multicols}
\caption{$2\times 2 \times 2$ tensors and their flattened frontal slices}
\end{center}
\end{figure*}
\vspace{-0.2in}
$$A=\left ( \begin{array}{ccccc}
a_{111} & a_{121} & \vdots & a_{112} & a_{122}\\
a_{211} & a_{221} & \vdots & a_{212} & a_{222}\end{array} \right ), \quad
B=\left ( \begin{array}{ccccc}
1 & 0 & \vdots & 0 & 1\\
0 & 1 & \vdots & 1 & 0\end{array} \right ).$$

For example, take  $A=(a_{ijk}),$ $ i, j, k=1, 2$. Then
$$\begin{array}{lll}
\1per (A) & = & a_{ 111}a_{ 222} +a_{211 }a_{122 }+a_{121 }a_{212 }+
a_{ 221}a_{112 },\\
\vspace{-.1in}
& & \\
\2per (A)& = & a_{ 111}a_{ 221} a_{122 }a_{212}+a_{211 }a_{121 }a_{112 }
a_{ 222 }.
\end{array}$$

Note that we use the same symbol for the tensor and its frontal slice flattening as there is no confusion will be caused in this paper.

By a (0,1)-tensor, we mean a tensor in which every entry is either 0 or 1. For the $2\times 2\times 2$ (0,1)-tensor $B$  in Fig.\,1 on the right hand side, $\1per (B)=0$, $\2per (B)=1$. If $I_3$ is the $3\times 3\times 3$ identity tensor, i.e.,
the entry in the $(i,i,i)$ position is 1 for every  $i$ and everywhere else is 0, then
$\1per (I_3)=1$ and $\2per  (I_3)=0$.
Let $J_3$ be the $3\times 3\times 3$ tensor of 1s (i.e., all entries are 1). Then
$\1per (J_3)=9\cdot 4 =(3!)^2=36$ and $\2per  (J_3)=6\cdot 2=12$.

1-permanent and 2-permanent are the most important permanent functions of tensors.
We write $\per (\cdot)$ for $\1per (\cdot $) and $\Per (\cdot )$ for $\2per (\cdot $). We simply call 1-permanent {\bf \em permanent} in the sense of (\ref{Def:1-per})  unless otherwise stated.

The following observations are immediate for permanents (i.e., 1-permanents):
(i) the permanent function of  tensors is linear with respect to each hyperplane; (ii) interchange of two hyperplanes does not change the permanent; (iii) if
 $A=(a_{i_1i_2\dots i_d})$ is a tensor of order $d$ and dimension $n$ and $A^{\sigma}=(a_{i_{\sigma (1)}i_{\sigma (2)}\dots i_{\sigma (d)}})$ is the $\sigma$-transpose of $A$, where  $\sigma \in S_d$, then
$A$ and $A^{\sigma}$ have the same permanent; (iv)
if  $A=(a_{i_1i_2\dots i_d})$ and $B=(b_{i_1i_2\dots i_d})$ are nonnegative tensors of the same size and  $A\leq B$ entrywise, then the permanent of $A$ is less than or equal to the permanent of $B$; and
  (v) the Laplace expansion theorem  holds.

The classic Frobenius-K\"{o}nig theorem (see, e.g., \cite[p.\,158]{ZFZbook11}) states that for an $n\times n$ nonnegative matrix $A$, the permanent of $A$ vanishes if and only if $A$ contains an $r\times s$ zero submatrix such that
$r+s=n+1$. 
The following result is an analog for tensors.

\begin{pro}[Dow and Gibson \cite{DowGib87}]
Let
$A$ be an
$n_1\times  n_2\times  \dots \times
n_d$ tensor.
If $A$ contains an $m_1\times m_2\times \cdots \times m_d$ zero sub-tensor such that $\sum_{k=1}^dm_k=1+\sum_{k=2}^dn_k$, then $\per (A)=0$; but not conversely.
\end{pro}

Consequently, for a nonnegative tensor of order $d$ and dimension $n$, if 1-permanent is positive, then
 $\sum_{k=1}^d m_k\leq (d-1)n.$
 The positivity of the permanent of a nonnegative tensor of order $d$ and dimension $n$ is characterized in terms of term rank in \cite{BWSZ15}: it is positive if and only if the term rank is $n$ .

Lower and upper bounds for the permanents of (0,1)-tensors (or matrices, or Latin squares, etc) are always interesting and challenging. Shown below is a lower bound of the permanent, given the number of zero entries.

\begin{pro}[Dow and Gibson \cite{DowGib87}] 
Let
$A$ be an
$n\times  n\times  \dots \times
n$ (0,1)-tensor of order $d$ with exactly $t$ 0s. Then
$\per (A)\geq (n^{d-1}-t)\left ( (n-1)!\right )^{d-1}.$
\end{pro}

The well-known Minc-Bregman theorem on a (0,1)-matrix  gives an upper bound for the permanent of the (0,1)-matrix in terms of the numbers of 1s on each row (or column). For tensors, we have the following.

\begin{pro}[Dow and Gibson \cite{DowGibAMS87}] Let
$A=(a_{ijk})$ be an
$n\times  n\times
n$ (0,1)-tensor. Let $r_i=\sum_{j,k}a_{ijk}$ for $i=1, 2, \dots, n$. Then
the Minc-Bregman type inequality for 1-permanent holds:
$$ {\per} (A)\leq \prod_{i=1}^n(r_i!)^{1/r_i}.$$
\end{pro}

\begin{pro}[Dow and Gibson \cite{DowGibAMS87}]
Let
$A=(a_{ijk})$ be an
$n\times  n\times
n$ (0,1)-tensor. Let $r_{ij}=\sum_{k}a_{ijk}$ for $i, j=1, 2, \dots, n$.  Then
the Minc-Bregman type inequality for 2-permanent holds:
$$ \Per (A)\leq \prod_{i=1}^n(r_{ij}!)^{1/r_{ij}}.$$
\end{pro}

A {\bf \em permutation tensor}  is a (0,1)-tensor in which every hyperplane contains one and only one 1. In particular, the usual permutation matrices are permutation tensors of order $2$; the identity tensor $I_n$ (all entries on the main diagonal $(i,i,i)$, $i=1, 2, \dots, n$,   are 1)  is a permutation tensor of order 3. (Note: permutation tensors are defined differently in the literature; see, e.g., \cite{LL14Per}.)
 Let $\Omega_n^d$ be the convex hull of the permutation tensors of order $d$ and dimension $n$.
An analog of the Van der Waerden conjecture (see, e.g.,  \cite{Minc83})  for tensors is surely appealing.
Dow and Gibson \cite{DowGib87} conjectured that
if
$A=(a_{i_1i_2\cdots i_d})\in \Omega_n^d$, then
$${\per}(A)\geq (n!/n^n)^{d-1}$$ with equality if and only if $A=(1/n^{d-1}) J_n,$
where  $J_n$  is the tensor of all 1s.

This is disproved by Taranenko \cite[Proposition 4, p.\,590]{Tar16}.
Taranenko presented as many as 13 conjectures concerning permanents and stochastic polytopes in \cite{Tar16}. We single out a couple that are easily stated and understood.

\begin{con}[Taranenko  \cite{Tar16}]
Let $A$ be an $\overbrace{n\times   \cdots \times n}^d$ line-stochastic  tensor. If  $d$ is even, then $\per (A)>0$.
\end{con}

\begin{con}[Taranenko \cite{Tar16}] 
Let $A$ be an $\overbrace{n\times   \cdots \times n}^d$ line-stochastic  tensor. If  $n$ is odd, then $\per (A)>0$.
\end{con}

For more discussions on this, see Theorems 19 and 22 of \cite{Tar16}.

\subsection{$k$-permanent and  the Hadamard product}

Let $A=(a_{i_1i_2\dots i_d})$ be a tensor of order $d$. For $1\leq k\leq d$, let $f=d-k$. If we fix
$f$ of the indices $i_1, i_2, \dots,  i_d$ and let the rest $k$  indices vary, then we obtain a sub-tensor of $A$. We call such a sub-tensor a {\em $k$-plane}  of $A$ (see \cite{DowGib87, Tar16}). 1-plane (1 free index) is referred to as a {\bf \em line} (or {\em  fiber} or {\em tube}); 2-plane (2 free indices) is simply a {\bf \em plane}; a $(d-1)$-plane  of an  order $d$ tenor  is  usually called a {\bf \em hyperplane}.

Dow and Gibson  \cite{DowGib87} defined the {\bf \em $k$-permanent of $A$}, denoted by $k$-$\per (A)$,  to be
the sum of all possible products of $n^k$ entries of $A$ so that no two entries are taken from the  same $(d-k)$-plane \cite[p.\,142]{DowGib87}. If such a selection of entries of $A$ does not exist, then we write $\kper (A)=0$.

\begin{rem}\label{Rem6}
{\rm
The existence of such selections of the entries of $A$ is
extensively studied (see, e.g., \cite{BrCs91, CF93, JurRys68, MarTar, Sch78}) and it is in the area of configurations and block designs in combinatorics (see, e.g., \cite{CD06}).}
\end{rem}

For $2\times 2\times 2$ tensors, we have demonstrated $\1per$ and $\2per$ in the previous examples. Permanents defined by (\ref{Def:1-per}) always exist. Let $A=(a_{ijst})$ be a $2\times 2\times 2\times 2$ tensor.  The 2-per($A$) is the sum of products of  $n^k=2^2=4$ entries of $A$ that are not in the same $d-k=4-2=2$-plane. Such a selection of entries is impossible for four sequences $i, j, s, t$ of length 4 whose components are 1 or 2:
$a_{i_1i_2i_3i_4} a_{j_1j_2j_3j_4} a_{s_1s_2s_3s_4} a_{t_1t_2t_3t_4}$. Thus,
$$\mbox{2-per ($A$)} =  \sum a_{i}a_{j}a_{s}a_{t}  =0.$$

 Let $A=(a_{i_1i_2\dots i_d})$ be an $\overbrace{n\times   \cdots \times n}^d$ tensor and let $1\leq k<d$.  A {\em $k$-per diagonal} of $A$ consists of $n^k$ entries of $A$; each  entry is from a $(d-k)$-plane and no two entries are from the same $(d-k)$-plane.
  A {$1$-per diagonal} is simply called a {\bf \em diagonal}; that is, a diagonal of a tensor of dimension $n$  consists of $n$ entries, no two are from the same hyperplane.
   For $d=2$, $k=1$, a 1-per diagonal of $A$ consists of $n$ entries of $A$ from different lines (i.e., rows and columns). For $d=3$, $k=1$, a 1-per diagonal of $A$ consists of $n$ entries of $A$, each of which is from 1-plane, no two fall on the same plane.
 For $d=3$, $k=2$, a 2-per diagonal of $A$ consists of $n^2$ entries of $A$ each plane contains exactly  $n$  entries of $A$.

We say that $A$ is {\em $k$-per feasible} if it is possible to choose
{$n^k$} entries of $A$, no two in the  same $(d-k)$-plane. Such a selection of the $n^k$ entries comprises of a { $k$-per diagonal} of $A$.  The $k$-per diagonal of $A$ can be extracted by the Hadamard (Schur or entrywise)  product of $A$ with a (0,1)-tensor $D$ of the same size (order and dimension) as $A$ in which the $k$-per diagonal entries of $D$ in the same positions as the $k$-per diagonal of $A$ are 1s and 0s elsewhere.
We call such a (0,1)-tensor $D$ a {\em  $k$-per  index tensor}
(or {\em a $k$-per cell}). That is, a $k$-per index tensor
is a (0,1)-tensor of order $d$ and dimension $n$ which contains $n^k$ 1s
so that no two 1s are located in the same $(d-k)$-plane.
Let $\prod(A\circ D)$
be the product of the $k$-per diagonal entries of $A$ indexed by $D$. Denote by $\mathcal{P}_{d,n,k}$, or simply $\mathcal{P}_{k}$, the set of $k$-per index tensors.
(Note: again, the existence of a $k$-per index tensor for a given $k$  is a problem of configuration which is not a concern of this paper.  For the study  of the existence of (0,1)-tensors with a fixed number of 1s on a $k$-plane, see, e.g., \cite{CF93, CD06, Sch78}.)

  We formulate the $k$-permanent of tensor $A$ \cite{DowGib87} as follows.

\begin{pro} Let $A=(a_{i_1i_2\dots i_d})$ be an $\overbrace{n\times   \cdots \times n}^d$ tensor, $1\leq k<d$.  Then
$$\kper (A)=\sum_{D\in  \mathcal{P}_{k}} \prod (A\circ D).$$
\end{pro}

\begin{pro} Let $A=(a_{i_1i_2\dots i_d})$ be an $\overbrace{n\times   \cdots \times n}^d$ tensor, $1\leq k<d$.  Then
$$\kper (cA)={c^{n^k}} \left (\kper (A)\right ), \;\; \mbox{where $c$ is a constant}.$$
\end{pro}

The following result states that every $k$-per can be converted to a 1-per.

\begin{pro}[Dow and Gibson \cite{DowGib87}]
 Let $A$ be a tensor of order $d$ and dimension $n$, $1\leq k<d$.  Then there exists an $n^k\times n^k \times \cdots \times n^k$ tensor $B$ of  order ${d\choose k}$ whose nonzero entries are equal to the nonzero entries of  $A$  such that
$$\kper (A)=\1per (B). $$
\end{pro}

\begin{rem}{\rm
 Different generalizations of the  permanents from  matrices  to tensors exist. Taranenko \cite{Tar16} defined $r$-permanents, ${\rm per}_r$, of multidimensional matrices  by the Maximum Distance Separable (MDS) codes with distance $r$. In \cite{Tar16}, the permanent  is in fact the $d$-permanent, that is,
$\per (A)={\rm per}_d(A)$, which is the same as the 1-permanent in \cite{DowGib87}, namely
our (\ref{Def:1-per}).
More generally, if $r+s = d+1$,
then the $r$-permanent in  \cite{Tar16}
is just the $s$-permanent in \cite{DowGib87}.
For $n\times n\times n$ tensors, the 2-permanents defined in \cite{DowGib87} and in
\cite{Tar16} turn out to be the same, namely, 2-per ($A$)= ${\rm per}_2(A)$.
However, for $2\times 2\times 2\times 2$ tensors,
$\2per (A) \not = {\rm per}_2(A)$.}
\end{rem}

\subsection{Permanent and the Hamming distance}


Let $x=(x_1, x_2, \dots, x_n)$, $y=(y_1, y_2, \dots, y_n)\in \Bbb R^n$.
The Hamming distance, denoted by ${\rho}(x, y),$ of $x$ and $y$ is the number of nonzero components of $x-y$. Take
$x=(1, 2, 3),$ $ y=(1, 3, 2)$. Then  $x-y=(0, -1, 1)$. Thus, $\rho(x,y)=2$.

Denote
$\mathcal{I}_n^d=\{(i_1, i_2, \dots, i_d)\}$, where $1\leq i_k\leq n$ for $k=1, 2, \dots, d$.
Let $A=(a_{i_1i_2\dots i_d})$. We write  $A=(a_{i})$, where $a_{i}=a_{i_1i_2\cdots i_d}$,
  $i=(i_1, \dots, i_d)\in \mathcal{I}_n^d$.
  For $i=(i_1, \dots, i_d)$ and $j=(j_1, \dots, j_d)$,  $\rho(i, j)=d$ implies
  that the corresponding components of $i$ and $j$ are pairwise  distinct.
  Taranenko \cite{Tar16} defined the permanent of a tensor using Hamming distance, which is essentially the same as (\ref{Def:1-per}), i.e., the 1-permanent of \cite{DowGib87}, namely  the
  $d$-permanent of \cite{Tar16}.

\begin{pro} Let $A=(a_{i_1i_2\dots i_d})$ be an $\overbrace{n\times   \cdots \times n}^d$ tensor.  Then
\begin{equation}\label{pro2}
\per (A)=\sum_{\substack{\a^{1},\, \a^{2}, \dots, \a^{d}\in\mathcal{I}_n^d \\
 \rho (\a^{i},\, \a^{j})=d,\; i\not = j}} a_{\a^{1}}a_{\a^{2}}\cdots a_{\a^{n}}.
 \end{equation}
\end{pro}

Let $A=(a_{ijst})$ be a $2\times 2\times 2\times 2$ tensor. Then $\per (A)$  is the sum of all products of  $n^k=2^1=2$ entries of $A$ that are not in the same $d-k=4-1=3$-plane.
It follows from (\ref{pro2})  that
\begin{eqnarray*}
\per (A)& = & \sum_{\rho(\a, \b)=4} a_{\a}a_{\b}  \\
  & = & a_{1111}a_{2222}+a_{1112}a_{2221}+a_{1121 }a_{2212 } + a_{1211 }a_{2122 }\\
  & & +a_{2111}a_{1222} +a_{1122 }a_{2211 } +
  a_{ 1212}a_{2121 }+a_{ 1221}a_{2112 }.
\end{eqnarray*}

Note that the  $\2per (A)$  is the sum of all products of  $n^k=2^2=4$ entries of $A$ that are not in the same $d-k=4-2=2$-plane. It is impossible for four sequences $\a, \b, \gamma, \delta$ of length 4 whose components are 1 or 2 to have $\rho(p, q)\geq 3$ for all pairs $p$ and $q$ from $\{\a, \b, \gamma, \delta\}$. Thus, $\Per (A)=0$.

\subsection{The permanents of 3D matrices (i.e., $n\times n\times n$ tensors)}

For $\a, \b \in S_n$, we may regard $\a$ and $\b$ as
sequences in $\Bbb R^n$: $\a=(\a(1), \dots, \a(n))$ and $\b=(\b(1), \dots,
\b(n))$.  If $\rho (\a, \b)=n,$ then $\a(i)\neq \b (i)$ for all $i$.

Let $A=(a_{ijk})$ be a 3D matrix (i.e., a tensor of order 3 and dimension $n$), namely, $A$ is an ${n\times n\times n}$ tensor. A {\bf \em  diagonal} of $A$ consists of $n$ entries, each of which is taken from a plane and no two entries are from the same plane (as $d-k=3-1=2$);
A {\bf \em triagonal} (\cite[p.\,181]{FiSw85}) of $A$  consists of $n^2$ entries, each of which is taken from a line  and no two   are from the same line. 
Then $\1per (A)$  is
the sum of products of diagonal entries. Thus \cite{{DowGibAMS87}},
\begin{equation}
\1per (A)=\sum_{\a,\, \b\in S_n} \prod_{i=1}^n a_{i\a(i)\b(i)},
\end{equation}
while $\2per (A)$  is
the sum of products of triagonal  entries. So,
$$\2per (A) =\sum \prod (\mbox{$n^2$ entries of $A$; no two are colinear}).  $$

For  $3\times 3\times 3$ tensors, 1-permanent is the sum of all products of 3 elements, no two are on the same plane  (frontal, lateral or horizontal \cite{KB09}), i.e., one element from each plane, while  2-permanent  is the sum of all products of
9 entries any two of which are non-collinear (in any direction), i.e., one entry from each line.

Let $A=(a_{ijk})$ be an ${n\times n\times n}$ tensor. If we
 denote (or label) the $k$th frontal page of $A$ by  $\pi_k\in S_n$, we can write
a {triagonal}  ${a}_{\pi}$  of $A$ as
$${a}_{\pi}=(a_{\pi_1}, a_{\pi_2}, \dots, a_{\pi_n})$$
where $\rho (\pi_i, \pi_j)=n$ whenever $i\neq j$. Let
$\mathcal{D}({a}_{\pi})=\prod_{i=1}^n \mathcal{D}({a_{\pi_i}})$ for the product of the triagonal entries.    Then (see, e.g.,  \cite{ChangPZ16, CLN14}), we have

\begin{pro} Let $A=(a_{ijk})$ be an ${n\times n\times n}$ tensor. Then
$$\Per (A)=\mbox{\rm 2-per ($A$)}=\sum_{\pi} \prod \mathcal{D}({a}_{\pi})
=\sum_{\substack{\pi_1, \pi_2, \dots, \pi_d\in S_n \\
\rho (\pi_i, \pi_j)=n, \,i\neq j}} \; \prod_{i=1}^n  \mathcal{D}(a_{\pi_i}).$$
\end{pro}

It is easy to see  that there are $n^2\cdot (n-1)^2 \cdots 2^2\cdot 1^2=(n!)^2$ permutation tensors  of order 3 and dimension $n$.
Let $L_n$ be the number of Latin squares of order $n$ and let
  $J_n^3$ denote  the $n\times n\times n$  tensor of all 1s. Then $\Per (J_n^3)$ is equal to the number of triagonals of $A=(a_{ijk})$.
Observe that  every triagonal of $J_n^3$ corresponds solely to a   Latin square of order $n$ (see, e.g.,
\cite{JurRys68}).

\begin{pro}
A 3D matrix of dimension $n$ has $L_n$ triagonals.
\end{pro}

\section{Stochastic tensors}

\subsection{Line, plane, $k$-stochastic, and permutation tensors}

Recall the celebrated Birkhoff-von Neumann theorem  on the polytope of doubly stochastic matrices
(see, e.g., \cite[p.159]{ZFZbook11}). It states that
an $n\times n$ matrix is doubly stochastic if and only if it is a convex combination of some $n\times n$ permutation matrices.
This result is about the matrices that are 2-way stochastic. What would be the mathematical objects that are
3-way stochastic?

 Let $A=(a_{ijk})$ be an $n\times n\times n$ tensor.  $A$ is said to be {\em triply line stochastic} \cite{CLN14} (or
 {\em   stochastic semi-magic cube} \cite{Mayathesis04}) if
all $a_{ijk}\geq 0$ and
$$\sum_{i=1}^n a_{ijk}=1, \quad \sum_{j=1}^n a_{ijk}=1, \quad  \sum_{k=1}^n a_{ijk}=1.$$
That is,   each of horizontal, lateral and frontal slices (see \cite{KB09})  is a doubly stochastic matrix. For a nonnegative tensor $A=(a_{i_1i_2\dots i_d})$
of order $d$ and dimension $n$, we say $A$ is {\em line-stochastic} \cite{FiSw85}
  if the sum of  the entries on each line is   1:
$$\sum_{i=1}^n a_{\cdots i\cdots }=1.$$
Equivalently, every plane (i.e., 2-plane) of $A$ is   doubly stochastic, namely, for $e=(1, 1, \dots, 1)^t \in \Bbb R^n$,
every $n\times n$ matrix with $(i, j)$ entry $a_{\cdots i\cdots j \cdots}$ satisfies
$$  (a_{\cdots i\cdots j \cdots})e=e, \quad  e^t (a_{\cdots i\cdots j \cdots}) =e^t.$$
We say that $A$ is {\em plane-stochastic} \cite{BrCs75laa} if the sum of all elements on every plane is equal to 1, that is,
$$\sum_{i, j=1}^n a_{\cdots i\cdots j \cdots}=1.$$

More generally, let $A$ be a nonnegative tensor of order $d$ and dimension $n$ and let $1\leq k\leq d$. If the sum of the entries of $A$ on every $k$-plane is 1, then $A$ is said to be {\bf \em $k$-stochastic} (see, e.g.,  \cite{BrCs91, Sch78}). A $k$-stochastic (0,1)-tensor is called a {\em $k$-permutation tensor} (or a permutation tensor of degree $k$; for its existence, see
Remark \ref{Rem6}).
 Being line stochastic is 1-stochastic; being 2-stochastic is plane-stochastic; and a 1-permutation tensor is nothing but a {\bf \em line-permutation} tensor, while a 2-permutation tensor is a {\bf \em plane-permutation} tensor. In case of
$n\times n\times n$,  1-permutation tensor has 1s on its diagonal
and  2-permutation tensor has 1s on its triagonal. The $(d-1)$-permutation tensors (of order $d$ and dimension $n$)  are simply called {\bf \em permutation tensors} \cite{DowGib87}.

Let $P$ and $Q$ be $n\times n$ permutation matrices. We say that $P$ and $Q$ are
{\em diagonally disjoint} (or {\em Hadamard orthogonal})  if no 1 appears in the same (overlapping) position of $P$ and $Q$, that is, the Hadamard product
$P\circ Q=0$.  

\begin{pro}
Let $P_1, P_2, \dots, P_n$ be $n\times n$ permutation matrices and $\pi_1, \pi_2,$ $ \dots,$ $ \pi_n$ be  the
 corresponding elements (via group isomorphism) in  the symmetric group $S_n$. The following statements are equivalent:
\begin{enumerate}
\item The tensor with frontal slice flattening  $[P_1 | P_2 | \cdots | P_n]$ is an $n\times n\times n$ line-permutation tensor.
\item $P_1, P_2, \dots, P_n$ are mutually diagonally disjoint.
\item
$\rho(\pi_i,\pi_j)=n$ for all $i\not = j$.
\item   $P_1+P_2+\cdots + P_n=J$ (where $J$ is the  matrix of 1s).
\end{enumerate}
\end{pro}

For an analog for $n\times n\times n$ plane-permutation tensors,
let $Q_1, Q_2, \dots, Q_n$ be $n\times n$ permutation matrices. Then
the $n\times n\times n$ (0,1)-tensor $R$ with frontal slice flattening $[Q_1| Q_2|\cdots | Q_n]$ is a plane-permutation tensor if and only if each  of the
{\em plus-projections} (by adding the elements) $f_i(R)$,  $f_j(R)$, and $f_k(R)$  of $R$ along $i$, $j$, and  $k$-axes is an $n\times n$ permutation matrix.

\subsection{Polytopes of stochastic tensors}
The Birkhoff-von Neumann Theorem asserts that the set of the doubly stochastic matrices and the convex hull of the permutation matrices coincide. In other words,  the permutation matrices are precisely  the vertices (extreme points) of
the polytope of doubly stochastic matrices.
This is usually proven by the Frobenius-K\"{o}nig theorem (see, e.g., \cite[p.158]{ZFZbook11}).

The Birkhoff-von Neumann Theorem does not generalize to tensors of higher dimensions. The
$3\times 3 \times 3$ line-stochastic tensor $D$ in Fig. 2 is not a combination of
line-permutation tensors; in fact, it is an extreme point of the polytope of
$3\times 3 \times 3$ line-stochastic tensors. 
Moreover,
$\Per (D)=0$. Let

\begin{itemize}
\item
$\Delta_{n}^{\ell}$ be the convex hull of  $n\times n \times n$ line-permutation tensors.
\item
$\Delta_{n}^{\wp}$ be the convex hull of  $n\times n \times n$ plane-permutation tensors.

\item
$\Omega_{n}^{\ell}$ be the set of all $n\times n \times n$ line-stochastic  tensors.
\item
$\Omega_{n}^{\wp}$ be the set of all $n\times n \times n$ plane-stochastic  tensors.
\end{itemize}

The $\Delta$s and $\Omega$s are polytopes in  $\Bbb R^{n^3}$.
It is tempting  to know the structures and the numbers of the extreme points of
the polytopes $\Delta$s. 
Obviously,
$$\Delta_{n}^{\imath}\subseteq \Omega_{n}^{\imath}, \quad \mbox{where $\imath=\ell$ or $\wp$}.$$

The following facts are known or easy to obtain:

\begin{enumerate}
\item For $n=2$, $\Delta_{2}^{\ell}= \Omega_{2}^{\ell}$. That is to say, every $2\times 2\times 2$ line-stochastic tensor is a convex combination of the  two (0,1) line-stochastic tensors.

\item For $n=2$, $\Delta_{2}^{\wp}$ is a proper subset of $\Omega_{2}^{\wp}$.
Take $C=(c_{ijk})$ with $c_{211}=c_{121}=c_{112}=c_{222}=\frac12$, and 0 everywhere else.
$C$ is not a convex combination of the plane-permutation tensors.
 $\Omega_{2}^{\wp}$ has 6 extreme points, 4 of which  are (0,1)-tensors and 2 are non-(0,1) (with entries $1/2$); see  \cite{BrCs75laa, Che17, LKLC, Sch78}. The cube on the left in Fig.\,2 represents the tensor $C$, while shown below is its frontal slice flattening. (Likewise, the other cube  is for tensor $D$.)
$$C=\frac12  \left ( \begin{array}{ccccc}
0 & 1 & \vdots & 1 & 0 \\
1 & 0 &\vdots & 0 &1 \end{array} \right ).$$

\begin{figure*}[h] 
\begin{center}
\begin{multicols}{2}
    \includegraphics[width=1.4in,totalheight=1.2in]{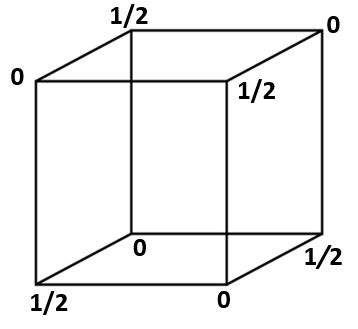}
    \includegraphics[width=1.5in,totalheight=1.2in]{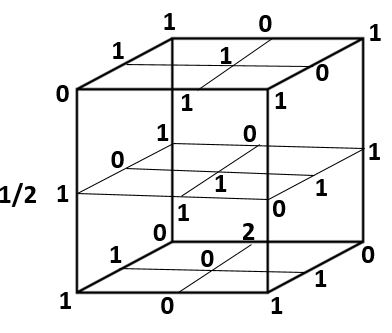}\par
    \end{multicols}
\caption{$C\in \Omega_{2}^{\wp}\setminus \Delta_{2}^{\wp}$ and  $D \in \Omega_{3}^{\ell}\setminus \Delta_{3}^{\ell}$}
\end{center}
\end{figure*}

\item For $n=3$,
the polytope $\Omega_{3}^{\ell}$ has  66 vertices, of which 12 are line-permutation tensors (due to the fact that there are 12 Latin squares of order 3), 54 are non-(0,1) (with entries $1/2$).
$\Delta_{3}^{\ell}$ is a proper subset of $\Omega_{3}^{\ell}$ because tensor $D$ is not a convex combination of line-permutation tensors (see, e.g., \cite{ChangPZ16}).

Moreover, for the  line-stochastic tensor $D$, we have $\Per (D)=0$.  This says, unlike the permanent of a doubly stochastic matrix, that  the 2-permanent (i.e., $\Per$)  of a triply line-stochastic tensor may vanish.
$$D=\frac12 \left  (
\begin{array}{ccccccccccc}
0 & 1 & 1 & \vdots  & 1 & 1 & 0 & \vdots &  1 & 0 & 1 \\
1 & 1 & 0 & \vdots  & 0 & 1 & 1 & \vdots &  1 & 0 & 1 \\
1 & 0&  1 & \vdots  & 1 & 0 & 1 & \vdots &  0 & 2 & 0
\end{array} \right ).$$

\item For $n=3$, $\Delta_{3}^{\wp}$ is a proper subset of $\Omega_{3}^{\wp}$.
 A complete list of the extreme points of $\Omega_{3}^{\wp}$, up to equivalence, is available in \cite{BrCs75laa}.

\end{enumerate}

\begin{que}{\rm
What would be the minimums  of the permanents on the sets $\Delta$s?}
\end{que}

\begin{pro}
Let $A=(a_{ijk})$ be an $n\times n\times n$ nonnegative tensor. If a plus-projection (by adding the elements) $f_i(A)$, $f_j(A)$, or $f_k(A)$  of $A$ along $i$, $j$, or $k$-axis
contains a $0$, then
$\Per (A)=0$. (The converse is not true.)
\end{pro}

If $R$ is a nonnegative tensor such that ($k$-per or) $\Per(R)>0$, then for any nonnegative tensor $S$ of the same size, ($k$-per, resp.) $\Per (R+S)\geq \Per (R)>0.$

\begin{pro} Let $P_0=\{A\in \Omega_{n}^{\ell} \mid {\Per} (A)=0\}$ and let $P$ and $ Q$ be in $ P_0$. Then either everything between $P$ and $Q$ is contained in $P_0$ (i.e., $tP+(1-t)Q\in P_0$ for all $0<t<1$), or nothing between $P$ and $Q$ lies in $P_0$ (i.e., $tP+(1-t)Q\not \in P_0$ for all $0<t<1$).
\end{pro}

\section{Generalized tensor functions}

Let $ A=(a_{ij})$ be an ${n\times n}$ matrix. Let $G$ be a subgroup of $S_n$ and $\chi$ be a character on $G$. The classic determinant, permanent, and generalized matrix function of $A$ are respectively defined  by
$$\det A=\sum_{\b \in S_n} (-1)^{{\rm {\rm sgn}} (\b )} \prod_{i=1}^n a_{{i} \b(i)}=\frac{1}{n!} \sum_{\a, \b \in S_n} (-1)^{{\rm sgn} (\a) {\rm sgn} (\b)} \prod_{i=1}^n a_{\a (i) \b(i)},$$
$$\per A=\sum_{\b \in S_n}  \prod_{i=1}^n a_{{i} \b(i)}=\frac{1}{n!}\sum_{\a, \b \in S_n}  \prod_{i=1}^n a_{\a (i) \b(i)},$$
$$d_{G}^{\chi}  A=\sum_{\b \in G} \chi (\b)  \prod_{i=1}^n a_{i \b(i)} =\frac{1}{|G|}\sum_{\a, \b \in G} \chi (\a)\chi (\b)  \prod_{i=1}^n a_{\a (i) \b(i)}.$$

For a tensor $A=(a_{i_1i_2\cdots i_d})$ of order $d$ and dimension $n$.
Cayley's  combinatorial (v.s. geometric)  hyperdeterminant 
of $A$ is defined to be 
\begin{equation}\label{CayleyDet}
\det A=\frac{1}{n!} \sum_{\pi_1, \dots, \pi_d\in S_n}{{\rm sgn}} (\pi_1) \dots {\rm sgn} (\pi_d) \prod_{i=1}^n a_{\pi_1(i)\cdots  \pi_d(i)}.
\end{equation}

The reader is referred to \cite{Gel94, Lim13, SSZ13, Hu2013} for hyperdeterminants or the determinants of multidimensional matrices (tensors).

For a tensor $A=(a_{i_1i_2\cdots i_d})$ of order $d$ and dimension $n$, the permanent (1-permanent) of $A$ is defined analogously as in (\ref{Def:1-per}) and (\ref{Def:1-perB}).
We now  give a try to extend the notation to generalized tensor functions.

 Let $A=(a_{i_1i_2\cdots i_d})$ be a tensor of order $d$ and dimension $n$. Let $G=(G_1, G_2, \dots, G_d)$ and $\chi =(\chi_1, \chi_2, \dots, \chi_d)$, where  $G_i$ is a subgroup of $S_n$ and $\chi_i$ is  a character on $G_i$, $i=1, 2, \dots, d$.
We define
\begin{equation}\label{CayleyGTF}
d_G^{\chi}(A)=\frac{1}{|G_1|}\sum_{\pi_1\in G_1,   \dots, \pi_d\in G_d} \;\chi_1(\pi_1)\cdots
\chi_d(\pi_d)\prod_{i=1}^n a_{\pi_1(i)\cdots  \pi_d(i)}.
\end{equation}

Apparently, the determinant (\ref{CayleyDet}) and permanent (\ref{Def:1-perB}) are special cases of (\ref{CayleyGTF}). Like 2-permanent for ${n\times n\times n}$ tensors,  we may define
2-$d_G^{\chi}$ as follows:
\begin{equation}\label{CayleyGMF}
\mbox{2-}d_G^{\chi}(A)=\sum_{\rho(\pi_i, \pi_j)=n, \,i\neq j} \; \prod_{i=1}^n {\chi_i}(\pi_i)a_{\pi_i}.
\end{equation}


Additionally, in regard to the $k$-permanent, we may define the $k$-generalized tensor functions ($k$-gtf). Let $f_k$ be a scalar-valued function defined on a domain that contains all $k$-per diagonals $A\circ D$ of $A$, where $D\in \mathcal{P}_k$ (see Sec.\,2.2). Then
\begin{equation}\label{gtf}
\mbox{$k$-gtf($A$)}=\sum_{D\in \mathcal{P}_k}f_k(A\circ D)\prod (A\circ D).
\end{equation}

The work of Merris \cite{Mer71} may be a hint for the study in this direction.

\bigskip

{\bf Acknowledgement.} The work  was done while the second author was visiting Shanghai University during his sabbatical leave from Nova Southeastern University. The work of Wang was partially supported  by the Natural Science Foundation of China (11571220); the work of Zhang was partially supported by an NSU PFRDG Research Scholar grant.  This expository  article was written based on the second author's presentation at ICMAA in Da Nang, Vietnam, June 14-18, 2017.
The authors  appreciate the communications with C. Bu, L. Cui, S. Hu, L. Qi, A. Taranenko, Y. Wei, and G. Yu during the preparation of the manuscript.

\end{document}